\documentclass[a4paper,10pt]{scrartcl}
\usepackage[english]{babel}
\usepackage[utf8]{inputenc}
\usepackage[T1]{fontenc}
\usepackage{lmodern} 
\usepackage{enumerate}
\usepackage[a4paper]{geometry}
\usepackage{amsfonts,amsmath,amssymb,amsthm,amsopn}
\usepackage{mathtools}
\usepackage{stmaryrd}
\usepackage{nicefrac}
\usepackage{exscale}

\usepackage[numbers,square]{natbib}
\usepackage{url} 
\usepackage[unicode]{hyperref}

\allowdisplaybreaks 

\newtheoremstyle{theo}
	{3pt} 
	{3pt} 
	{\itshape} 
	{} 
		{\bfseries} 
	{\\} 
	{ } 
	{\thmname{#1}\thmnumber{ #2.}\thmnote{ - #3}} 
\theoremstyle{theo}

\newtheorem{defi}{Definition}[section]
\newtheorem{lem}[defi]{Lemma}
\newtheorem{theo}[defi]{Theorem}
\newtheorem{cor}[defi]{Corollary}
\newtheorem{rem}[defi]{Remark}
\newtheorem{prop}[defi]{Proposition}

\newenvironment{bew}{\begin{proof}[\bfseries Proof:]}{\end{proof}}

\DeclareMathOperator{\bomega}{\overline{\Omega}}
\DeclareMathOperator{\romega}{\partial\Omega}
 
\DeclareMathOperator{\Lo}{L}
\DeclareMathOperator{\W}{W}

\DeclareMathOperator{\intd}{d\!}
\DeclareMathOperator{\divdot}{\!\cdot}

\newcommand{\Tm}{T_{max}}
\newcommand{\etd}{e^{t\Delta}}

\newcommand{\etsd}{e^{(t-s)\Delta}}
\newcommand{\etsdm}{e^{(t-s)(\Delta-1)}}
\newcommand{\intot}{\int\limits_0^{t}} 

\newcommand{\intomega}{\int\limits_{\Omega}\!} 

\newcommand{\sintomega}{\smallint_{\text{\tiny{$\Omega$}}}} 

\newcommand{\intromega}{\int\limits_{\romega}\!} 

\newcounter{gleichung}
\newcommand{\owncount}{\refstepcounter{gleichung}}

\author{Tobias Black\thanks{Institut f\"ur Mathematik, Universit\"at Paderborn, Warburger Str. 100, 33098 Paderborn, Germany; email: \mbox{tblack@math.upb.de}}}
\title{Boundedness in a Keller-Segel system with external signal production}

\begin{document}
\maketitle

\begin{abstract}
\noindent
{\textbf{Abstract: } We study the Neumann initial-boundary problem for the chemotaxis system
\begin{align*}
\left\{\begin{array}{c@{\,}l@{\quad}l@{\,}c}
u_{t}&=\Delta u-\nabla\divdot(u\nabla v),\ &x\in\Omega,& t>0,\\
v_{t}&=\Delta v-v+u+f(x,t),\ &x\in\Omega,& t>0,\\
\frac{\partial u}{\partial\nu}&=\frac{\partial v}{\partial\nu}=0,\ &x\in\romega,& t>0,\\
 u(x,0)&=u_{0}(x),\ v(x,0)=v_{0}(x),\ &x\in\Omega&
        \end{array}\right.
\end{align*}
in a smooth, bounded domain $\Omega\subset\mathbb{R}^n$ with $n\geq2$ and $f\in\Lo^\infty([0,\infty);\Lo^{\frac{n}{2}+\delta_0}(\Omega))\cap C^\alpha(\Omega\times(0,\infty))$ with some $\alpha>0$ and $\delta_0\in\left(0,1\right)$.

First we prove local existence of classical solutions for reasonably regular initial values. Afterwards we show that in the case of $n=2$ and $f$ being constant in time, requiring the nonnegative initial data $u_0$ to fulfill the property $\sintomega u_0\intd x<4\pi$ ensures that the solution is global and remains bounded uniformly in time. Thereby we extend the well known critical mass result by Nagai, Senba and Yoshida for the classical Keller-Segel model (coinciding with $f\equiv 0$ in the system above) to the case $f\not\equiv 0$.

Under certain smallness conditions imposed on the initial data and $f$ we furthermore show that for more general space dimension $n\geq2$ and $f$ not necessarily constant in time, the solutions are also global and remain bounded uniformly in time. Accordingly we extend a known result given by Winkler for the classical Keller-Segel system to the present situation.
}\\
{\textbf{Keywords:} }chemotaxis, boundedness, external signal production, critical mass\\
{\textbf{MSC (2010):} }35K55 (primary), 35A01, 35A02, 35Q92, 92C17
\end{abstract}
\pagebreak

\section{Introduction}\label{ch:c1}
%
In mathematical biology, PDE systems of the form
\begin{align}\label{KS}\tag{$KS$}
\left\{\begin{array}{c@{\,}l@{\quad}l@{\,}c}
u_{t}&=\Delta u-\nabla\divdot(u\nabla v),\ &x\in\Omega,& t>0,\\
v_{t}&=\Delta v-v+u,\ &x\in\Omega,& t>0,\\
\frac{\partial u}{\partial\nu}&=\frac{\partial v}{\partial\nu}=0,\ &x\in\romega,& t>0,\\
 u(x,0)&=u_{0}(x),\ v(x,0)=v_{0}(x),\ &x\in\Omega&
        \end{array}\right.
\end{align}
are widely being used to model the process of chemotaxis -- a biological phenomenon of oriented movement of cells in response to some kind of chemical substance. Systems of this type were introduced in 1970, when Keller and Segel proposed a mathematical model describing the aggregation of some types of bacteria (see \cite{KS70} and \cite{KS71}). The model given above is a special case of the system originally stated in the pioneering works. Therein $u(x,t)$ represents the density of the cells and $v(x,t)$ denotes the concentration of an attracting chemical substance at place $x$ and time $t$. The first equation of \eqref{KS} models the movement of the cells. This movement, while diffusive, also favors the direction toward higher concentration of the chemical substance. The second equation in \eqref{KS} models the assumptions that the chemical, while diffusing and degrading, is also consistently produced by the living cells themselves.

Similar variants to the system above have been used in modeling a wide array of biological phenomena, e.g. pattern formation in E. coli colonies (\cite{AMM10}) and cancer invasion of tissue (\cite{SRLC09}) to just name a few. For a broader variety and further impressions on the biological background of this type we refer to the survey articles \cite{BBWT15}, \cite{HP09} and \cite{Ho03}. Often the occurrence of self-organized patterns such as aggregation, is identified with the blow-up of the solution, i.e. the existence of some $T\in(0,\infty]$ such that $\limsup_{t\nearrow T}\|u\|_{\Lo^\infty(\Omega)}=\infty$. Accordingly, mathematical efforts are often focused on detecting unbounded solutions with finite time or infinite time blow up, or the lack of unbounded solutions. More generally one is interested in finding conditions on the initial data which either ensure or dismiss the existence of blow-up solutions all together.

The system \eqref{KS} has been thoroughly studied with regard to the boundedness of solutions. We briefly summarize some known results, where, if not stated otherwise, $\Omega\subset\mathbb{R}^n$ is an arbitrary, smooth and bounded domain and the initial values fulfill $u_0\in C^0(\bomega)$, $v_0\in C^1(\bomega)$ and are nonnegative. We denote by $(u,v)$ the corresponding maximally extended classical solution of \eqref{KS}:
\begin{itemize}
\item[If] $n=1:$ Then $(u,v)$ is global and bounded with regard to the $\Lo^\infty(\Omega)$-norm. (\cite{OY01})
\item[If] $n=2:$ If $\sintomega u_0\intd x<4\pi$ (or $8\pi$ in the radial symmetric setting), then $(u,v)$ is global and bounded with regard to the $\Lo^\infty(\Omega)$-norm. (\cite{NSY97}, \cite{gz98})\\
On the other hand, for any $m>4\pi$ with $m\not\in\lbrace 4k\pi\vert k\in\mathbb{N}\rbrace$ there exist initial data $u_0$,$v_0$ satisfying $\sintomega u_0\intd x=m$ such that the respective solution blows up in finite or infinite time. (\cite{hw01},\cite{ss01})

An additional important result is due to \cite{mizoguchi_winkler_13}, where it was shown -- in the radially symmetric setting -- that finite time blow up is a quite typical occurrence, in the sense that for each $p\in(0,1)$, $q\in(1,2)$ and a given $T>0$ the set of initial values leading to blow up before time $T$ is dense in $\{(u_0,v_0)\in C^0(\bomega)\times W^{1,\infty}(\Omega)\vert\ \text{radially symmetric and positive in }\bomega \}$ with respect the topology in $\Lo^p(\Omega)\times\W^{1,q}(\Omega)$.

\item[If] $n\geq3:$ It was proven in \cite{win10jde} that there exists a bound for $u_0$ in $\Lo^q(\Omega)$ and for $\nabla v_0$ in $\Lo^p(\Omega)$, with $q>\frac{n}{2}$ and $p>n$ such that the solution $(u,v)$ is global in time and bounded. This result has further been extended to the critical case $q=\frac{n}{2}$ and $p=n$. (\cite{cao2014global}) Regarding blow up of solutions it was shown in \cite{win10jde} -- in the radially symmetric setting -- that for any $m>0$ one can find initial data $u_0$ and $v_0$ satisfying $\sintomega u_0\intd x=m$, such that the solution blows up in finite or infinite time. Furthermore, it was proven in \cite{Win13pure} that -- similar to the two dimensional case -- blow up is also quite typical in space dimension three and above.
\end{itemize}
Similar results regarding the boundedness for the associated Cauchy problems on the whole space $\mathbb{R}^n$ have also been proven. See \cite{CP08} for instance, where dimensions $n\geq3$ are studied, or \cite{Miz13} for dimension $n=2$.

The main purpose of this work is to examine if corresponding statements hold when an \emph{additional external production of the signal chemical} is introduced to the system. More precisely, we shall study the initial-Neumann boundary value problem
\begin{align}\label{KSf}\tag{$KS_f$}
\left\{\begin{array}{c@{\,}l@{\quad}l@{\,}c}
u_{t}&=\Delta u-\nabla\divdot(u\nabla v),\ &x\in\Omega,& t>0,\\
\tau v_{t}&=\Delta v-v+u+f(x,t),\ &x\in\Omega,& t>0,\\
\frac{\partial u}{\partial\nu}&=\frac{\partial v}{\partial\nu}=0,\ &x\in\romega,& t>0,\\
 u(x,0)&=u_{0}(x),\ v(x,0)=v_{0}(x),\ &x\in\Omega&
        \end{array}\right.
\end{align}
in a bounded and smooth domain $\Omega\subset\mathbb{R}^n$ with $n\geq2,\tau>0,u_0\in C^0(\bomega), v_0\in\W^{1,\theta}(\Omega)$ for some $\theta>n$ and $f\in\Lo^\infty([0,\infty);\Lo^{\frac{n}{2}+\delta_0}(\Omega))\cap C^\alpha(\Omega\times(0,\infty))$ with some $\delta_0\in\left(0,1\right)$ and $\alpha>0$.

Let us briefly summarize the structure of this work along with our three main results:\vspace*{2mm}

\textbf{Local existence and uniqueness of classical solutions:}\\ Based on well known regularity estimates, we will first establish the local existence and uniqueness of classical solutions to \eqref{KSf}. 
\begin{theo}[Time local existence of classical solutions to \eqref{KSf}]\label{localexist}
Assume that $u_0\in C^0(\bomega)$ and $v_0\in\W^{1,\theta}(\Omega)$, for some $\theta>n$, are nonnegative. Furthermore suppose that $f\in\Lo^\infty([0,\infty);\Lo^{\frac{n}{2}+\delta_0}(\Omega))\cap C^\alpha(\Omega\times(0,\infty))$ for some $\alpha>0$ and $1>\delta_0>0$ is nonnegative. Then there exists a maximal existence time $\Tm\in(0,\infty]$ and a uniquely determined pair of nonnegative functions $(u,v)$ with the properties
\begin{align*}
u&\in C^0(\bomega\times[0,\Tm))\cap C^{2,1}(\bomega\times(0,\Tm))\\
v&\in C^0(\bomega\times[0,\Tm))\cap C^{2,1}(\bomega\times(0,\Tm))\cap \Lo_\text{loc}^\infty([0,\Tm);\W^{1,\theta}(\Omega)),
\end{align*}
solving \eqref{KSf} in the classical sense in $\Omega\times(0,\Tm)$. Moreover we obtain the alternative
\begin{align}\label{eq2.1}
\owncount
\text{either}\quad\Tm=\infty\ \text{ or}\quad\Vert u(\cdot,t)\Vert_{\Lo^\infty(\Omega)}+\Vert v(\cdot,t)\Vert_{\W^{1,\theta}(\Omega)}\rightarrow\infty\ \text{ as}\quad t\nearrow\Tm.
\end{align}
\end{theo}
The theorem above does not only give us local solutions to \eqref{KSf}, but also a criterion for their extensibility to global solutions. In view of the alternative \eqref{eq2.1}, the boundedness of $\|u\|_{\Lo^\infty(\Omega)}$ and $\|v\|_{\W^{1,\theta}(\Omega)}$ for all $t>0$ is sufficient in order to obtain global solutions.\vspace*{2mm}

\textbf{Boundedness of solutions in the two dimensional case:}\\
In this section we will limit ourselves to domains in $\mathbb{R}^2$. The reasoning behind this restriction is that in the case $n=2$ we are able to give an explicit value for the critical mass $m_u:=\sintomega u_0\intd x$, which -- in the two-dimensional case -- is the determining quantity in whether a solution is bounded globally in time or not. This behavior is well studied for Keller-Segel systems without external signal production -- see \cite{NSY97} and \cite{hw01} for instance. Following the methods of \cite{NSY97}, we show that the result concerning critical mass $m_u<4\pi$ for systems without external signal production -- mentioned above in the summary of known results -- extends to the case $f\not\equiv 0$ without a change in the critical mass, provided that $f$ is a nonnegative constant-in-time function of class $\Lo^{1+\delta_0}(\Omega)\cap C^\alpha(\Omega)$.
To be more precise we obtain the following:
\begin{theo}[Global existence for small initial mass]\label{globex2th}
Let $0<\alpha$, $0<\delta_0<1$ and $2<\theta<\frac{2+2\delta_0}{1-\delta_0}$. Furthermore assume that the functions $u_0\in C^0(\bomega)$, $v_0\in\W^{1,\theta}(\Omega)$ and $f\in\Lo^{1+\delta_0}(\Omega)\cap C^\alpha(\Omega)$ are nonnegative and $m_u:=\sintomega u_0\intd x<4\pi$. Then the unique classical solution $(u,v)$ to \eqref{KSf} is global in time and the quantities $\|u(\cdot,t)\|_{\Lo^{\infty}(\Omega)}$, $\|v(\cdot,t)\|_{\W^{1,\theta}(\Omega)}$ remain bounded uniformly in time.
\end{theo}

To prove this statement we define a generalized energy functional $$W(u,v)=\intomega\left[u\log u-uv+\frac{1}{2}(\vert\nabla v\vert^2+v^2)-fv\right]\intd x$$ and first show that this energy remains bounded for all times if $(u,v)$ is a classical solution to \eqref{KSf} and $m_u$ is strictly smaller than the threshold number. This bound will be the essential ingredient in the estimation process concerning the $\Lo^\infty(\Omega)$-norm of $u$ and the $\W^{1,\theta}(\Omega)$-norm of $v$, which by the alternative shown in Theorem \ref{localexist} will be sufficient to ensure global existence.\vspace*{2mm}

\textbf{Boundedness of solutions the higher dimensional case:}\\
The last section is dedicated to the more general case of $n\geq2$ and not necessarily time-independent $f\in\Lo^\infty([0,\infty);\Lo^{\frac{n}{2}+\delta_0}(\Omega))\cap C^\alpha(\Omega\times(0,\infty))$, for which we prove -- based on methods shown in \cite{win10jde} -- the boundedness of solutions if the three quantities $\|u_0\|_{\Lo^{\frac{n}{2}+\delta_0}(\Omega)},\|\nabla v_0\|_{\Lo^\theta(\Omega)}$ and $\|f\|_{\Lo^\infty([0,\infty);\,\,\Lo^{\frac{n}{2}+\delta_0}(\Omega))}$ are all bounded by some sufficiently small $\varepsilon>0$.
\begin{theo}[Global existence of small-data solutions]\label{globalexist3+}
Let $0<\delta_0<1$, $0<\alpha$, $n<\theta<\frac{n^2+2n\delta_0}{n-2\delta_0}$ and $1<r$. Then there exist constants $\varepsilon_0>0$ and $C>0$ with the following property:
If $u_0\in C^0(\bomega)$, $v_0\in\W^{1,\theta}(\Omega)$ and $f\in\Lo^\infty([0,\infty);\Lo^{\frac{n}{2}+\delta_0}(\Omega))\cap C^\alpha(\Omega\times(0,\infty))$ are nonnegative with
\begin{align}\label{globalexist3+eq1}
\owncount
\|u_0\|_{\Lo^{\frac{n}{2}+\delta_0}(\Omega)}\leq\varepsilon,\ \|\nabla v_0\|_{\Lo^\theta(\Omega)}\leq\varepsilon\text{ and }\|f\|_{\Lo^\infty([0,\infty);\,\Lo^{\frac{n}{2}+\delta_0}(\Omega))}\leq\varepsilon
\end{align}
for some $\varepsilon<\varepsilon_0$, then the classical solution $(u,v)$ of \eqref{KSf} exists globally with the quantities $\|u\|_{\Lo^\infty(\Omega)}$ and $\|v\|_{\W^{1,\theta}(\Omega)}$ remaining bounded for all times. Additionally these global small-data solutions remain small for all times in the sense that the solution satisfies
\begin{align}\label{globalexist3+eq2}
\owncount
\|u(\cdot,t)-\etd u_0\|_{\Lo^\infty(\Omega)}\leq C\varepsilon^2 e^{-\frac{\lambda_1}{r} t}+C\varepsilon^2\ \text{ for all }t>1,
\end{align}
as well as
\begin{align}
\owncount\label{globalexist3+eq3}
\left\Vert\nabla\!\left(v(\cdot,t)-e^{ \frac{t}{\tau}(\Delta-1)} v_0-\frac{1}{\tau}\intot \etsdm e^{s\Delta}u_0\intd s\right)\right\Vert_{\Lo^\theta(\Omega)}\hspace*{-10pt}\leq C\varepsilon^2 e^{-\frac{\lambda_1 t}{r}}+C\varepsilon\ \text{for all $t>1$.} 
\end{align}
Here $\lambda_1>0$ denotes the first nonzero eigenvalue of $-\Delta$ in $\Omega$ with regard to the homogeneous Neumann boundary conditions, and by $\left(\etd\right)_{t\geq0}$ we abbreviate the Neumann heat semigroup in $\Omega$.
\end{theo}
In the case $n=2$ the theorem above can be interpreted as an extension of our result from Theorem \ref{globex2th} in the sense that we now allow time-dependent external signal production $f(x,t)$, but no precise information on the explicit threshold number is gained in this case. 

For $n\geq3$ the theorem generalizes the statement from \cite{win10jde}, mentioned at the start of this section. Furthermore, \eqref{globalexist3+eq2} and \eqref{globalexist3+eq3} show that the solution $(u,v)$ almost, with a small offset depending on $\varepsilon$, asymptotically behaves like the solution $(u_H,v_H)$ of the two linear parabolic equations $\frac{\partial}{\partial t}u_{H}=\Delta u_{H}$ and $\frac{\partial}{\partial t} v_{H}=\Delta v_{H}-v_{H}+u_{H}$ under the same initial and boundary data.

\noindent
Accordingly, for reasonably well-behaved signal production most of the known results regarding the boundedness of solutions for the unexcited system \eqref{KS} can be transferred directly to the extended system \eqref{KSf}. That these statements may not hold for less regular signal production was shown in \cite{TW11}, where the parabolic-elliptic Keller-Segel system
\begin{align*}
\left\{\begin{array}{c@{\,}l@{\quad}l@{\,}c}
u_{t}&=\Delta u-\nabla\divdot(u\nabla v)\ &x\in\mathbb{R}^n,& t>0,\\
0 &=\Delta v+u+f(x)\ &x\in\mathbb{R}^n,& t>0,\\
u(x,0)&=u_{0}(x)\ &x\in\mathbb{R}^n,&
        \end{array}\right.
\end{align*}
with a Dirac distributed signal production $f(x)=f_0\delta(x)$ was considered. It was shown in a radially symmetric setting, that for any choice of $f_0>0$ certain generalized solutions, so called radial weak solutions, blow up immediately. This result suggests that the exponent $\delta_0=0$ is critical in \eqref{KSf} and needs to be discussed further for the higher dimensional case $n\geq3$.

\noindent
Throughout this work $n\geq2$ will always represent the space dimension, $\Omega\subset\mathbb{R}^n$ a bounded and smooth domain and furthermore the relations $\theta>n$, $\tau>0$ and $1>\delta_0>0$ will also always hold.
%
\setcounter{gleichung}{0} 
\section{Local existence and uniqueness}\label{cloc}
The proof of our local existence result is based on the method seen in \cite{win10ctax}, which relies heavily on the well known Neumann heat semigroup estimates. For a version suitable to our context we refer to \cite[Lemma 1.3]{win10jde}. We also require some results for fractional powers of sectorial operators shown in \cite{hen81}. For a more general introduction to sectorial operators and fractional powers see \cite{fr69} for instance.

\subsection{Local existence of classical solutions and uniqueness}\label{clocs2}
Before stating the main result of this section we give a very short proof of the fact that the system \eqref{KSf} conserves the mass with respect to $u$, which we will use multiple times later on.
\begin{lem}\label{conservmass}
If $(u,v)$ is a classical solution to \eqref{KSf} in $\Omega\times(0,T)$, then $\|u(\cdot,t)\|_{\Lo^1(\Omega)}=\|u_0\|_{\Lo^1(\Omega)}$ holds for all $t\in[0,T)$.
\end{lem} 
\begin{bew}
Integrating the first equation of \eqref{KSf} and applying the divergence theorem, we obtain
\begin{align*}
\frac{\intd}{\intd t}\intomega u\intd x&=\intomega\Delta u\intd x-\intomega\nabla\divdot(u\nabla v)\intd x=\intromega\frac{\partial u}{\partial \nu}\intd x-\intromega u\frac{\partial v}{\partial \nu}\intd x =0\,\text{ for }t\in(0,T),
\end{align*}
from which $\|u(\cdot,t)\|_{\Lo^1(\Omega)}=\text{const.}=\| u_0\|_{\Lo^1(\Omega)}$ follows.
\end{bew}

For the sake of notation let us briefly recall Young's inequality.
\begin{lem}[Young's inequality]\label{young}
Let $a,b,\varepsilon>0$ and $1<p,q<\infty$ with $\frac{1}{p}+\frac{1}{q}=1$. Then
\begin{align*}
ab\leq\varepsilon a^p+C(\varepsilon,p,q) b^q,
\end{align*}
where $C(\varepsilon,p,q)=(\varepsilon p)^{-\frac{q}{p}}q^{-1}$.
\end{lem}

\pagebreak

\begin{proof}[\textbf{Proof of Theorem \ref{localexist}:}]
First we observe that since $n<\frac{n^2+2n\delta_0}{n-2\delta_0}$ we may assume
\begin{align}\label{thetareq}
\owncount
n<\theta<\frac{n^2+2n\delta_0}{n-2\delta_0}
\end{align}
without loss of generality. Otherwise we take some $\theta^*<\theta$ fulfilling \eqref{thetareq}, apply the semigroup estimates below for $\theta^*$ instead of $\theta$ and make use of the embedding $\W^{1,q}(\Omega)\hookrightarrow\W^{1,p}(\Omega)$ for all $q\geq p$, to then estimate the $\W^{1,\theta^*}(\Omega)$-norm from above by the $\W^{1,\theta}(\Omega)$-norm.

\textbf{Existence:} To prove existence, we work along the lines of a basic fixed point argument. First we claim that for all $R>0$ there exists $T:=T(R)>0$ such that \eqref{KSf} has a mild solution in $\Omega\times(0,T)$ if we require $\Vert u_0\Vert_{\Lo^\infty(\Omega)}\leq R$ and $\Vert v_0\Vert_{\W^{1,\theta}(\Omega)}\leq R$ in addition to the assumptions above.
To substantiate this claim, recalling the Neumann heat semigroup estimates, we max fix $K>0$ such that $\Vert e^{\frac{t}{\tau}\Delta} w\Vert_{\W^{1,\theta}(\Omega)}\leq K\Vert w\Vert_{\W^{1,\theta}(\Omega)}$ holds for all $w\in\W^{1,\theta}(\Omega)$ and define, for sufficiently small $T\in(0,1)$ to be specified later, the Banach space
\begin{align*}
X:=C^0([0,T];C^0(\bomega))\times C^0([0,T];\W^{1,\theta}(\Omega))
\end{align*}
with its closed subset
\begin{align*}
S:=\left\lbrace(u,v)\in X\,\left\vert\,\Vert u\Vert_{\Lo^\infty((0,T);\Lo^\infty(\Omega))}\leq R+1,\Vert v\Vert_{\Lo^\infty((0,T);\W^{1,\theta}(\Omega))}\leq KR+1\right.\right\rbrace.
\end{align*}
For convenience of notation we will write $u(s),v(s)$ and $f(s)$ instead of $u(\cdot,s),v(\cdot,s)$ and $f(\cdot,s)$ in the following estimation processes. Introducing the map 
\begin{align*}
&\Phi(u,v)(t):=\begin{pmatrix}
\Phi_1(u,v)(t)\\
\Phi_2(u,v)(t)
\end{pmatrix}\\
:=&\begin{pmatrix}
\etd u_0-\intot\etsd\nabla\divdot(u(s)\nabla v(s))\intd s\\
e^{\frac{t}{\tau}(\Delta-1)}v_0+\frac{1}{\tau}\intot\etsdm u(s)\intd s+\frac{1}{\tau}\intot\etsdm f(s)\intd s
\end{pmatrix}
\end{align*}
for $(u,v)\in S$ and $t\in[0,T]$, we will show that, for sufficiently small $T$, $\Phi$ maps from $S$ onto itself. To this end we estimate both components of $\Phi$ separately, starting with $\Phi_1$. Obviously it is
\begin{align}\label{eq2.2}
\owncount
\Vert\Phi_1(u,v)(t)\Vert_{\Lo^\infty(\Omega)}\leq\Vert\etd u_0\Vert_{\Lo^\infty(\Omega)}+\intot\left\Vert\etsd\nabla\divdot(u(s)\nabla v(s))\right\Vert_{\Lo^\infty(\Omega)}\intd s.
\end{align}
Picking $p>\max\left\lbrace \frac{n\theta}{\theta-n},\frac{n^2+2n\delta_0}{n-2\delta_0}\right\rbrace$, we see that the interval $I:=(\frac{n}{p},\frac{1}{2}-\frac{n}{2}(\frac{1}{\theta}-\frac{1}{p}))$ is not empty. Choosing some $\alpha\in I$, we have $p\alpha>n$, and thus the fractional power $A^\alpha$ of the sectorial operator $A:=-\Delta+1$ with Neumann data in $\Lo^p(\Omega)$ satisfies $\Vert w\Vert_{\Lo^\infty(\Omega)}\leq C\Vert A^\alpha w\Vert_{\Lo^p(\Omega)}$ and $\Vert A^\alpha e^{\sigma\Delta} w\Vert_{\Lo^p(\Omega)}\leq C\sigma^{-\alpha}\Vert w\Vert_{\Lo^p(\Omega)}$ for all $w\in C_0^\infty(\Omega)$ and $\sigma>0$, with some suitable positive constant $C$ (cf. \cite[Theorems 1.6.1 and 1.4.3]{hen81}). This allows to estimate the second term of \eqref{eq2.2} according to
\begin{align*}
I_1&:=\intot\left\Vert\etsd\nabla\divdot(u(s)\nabla v(s))\right\Vert_{\Lo^\infty(\Omega)}\intd s\leq C_1\intot\left\Vert A^\alpha\etsd\nabla\divdot(u(s)\nabla v(s))\right\Vert_{\Lo^p(\Omega)}\intd s\nonumber\\
&\leq C_2\intot\left(\frac{t-s}{2}\right)^{-\alpha}\left\Vert e^{\frac{t-s}{2}\Delta}\nabla\divdot(u(s)\nabla v(s))\right\Vert_{\Lo^p(\Omega)}\intd s\quad \text{for all }t\in(0,T),\nonumber
\end{align*}
where here and below every $C_i$ denotes a suitable positive constant independent of $t$.
By the assumptions imposed on $\theta$ in \eqref{thetareq} and the choice of $p$ we have $p\geq\theta$. Thus we can utilize known semigroup estimates (\cite[Lemma 1.3]{win10jde}) and the fact $T<1$ to obtain
\begin{align}
I_1&\leq C_3\intot(t-s)^{-\alpha-\frac{1}{2}-\frac{n}{2}(\frac{1}{\theta}-\frac{1}{p})}\Vert u(s)\nabla v(s)\Vert_{\Lo^\theta(\Omega)}\intd s\nonumber\\&\leq C_3(R+1)(KR+1)\intot(t-s)^{-\frac{1}{2}-\alpha-\frac{n}{2}(\frac{1}{\theta}-\frac{1}{p})}\intd s\quad \text{for all }t\in(0,T)\nonumber.
\intertext{Herein, by the choice of $\alpha$, we have $-\frac{1}{2}-\alpha-\frac{n}{2}(\frac{1}{\theta}-\frac{1}{p})>-1$ and thus
}
\owncount\label{eq2.3}
I_1&\leq C_3(R+1)(KR+1)T^{\frac{1}{2}-\alpha-\frac{n}{2}(\frac{1}{\theta}-\frac{1}{p})}\quad \text{holds for all }t\in(0,T).
\end{align}
To treat the first term in \eqref{eq2.2} we can apply the maximum principle to easily deduce
\begin{align}\label{eq2.4}
\owncount
\Vert\etd u_0\Vert_{\Lo^\infty(\Omega)}\leq\Vert u_0\Vert_{\Lo^\infty(\Omega)}\leq R.
\end{align}
Now combining \eqref{eq2.2} -- \eqref{eq2.4} yields
\begin{align}\label{eq2.5}
\owncount
\Vert\Phi_1(u,v)(t)\Vert_{\Lo^\infty(\Omega)}\leq R+C_3(R+1)(KR+1)T^{\frac{1}{2}-\alpha-\frac{n}{2}(\frac{1}{\theta}-\frac{1}{p})}
\end{align}
for $t\in(0,T)$. To treat $\Phi_2$ we argue in a similar fashion. Since $1>\delta_0>0$ and $\theta>n\geq2$ we have $\frac{n}{2}+\delta_0=:q_0<\frac{n}{2}+1<\theta$. Additionally by \eqref{thetareq} we have $\theta<\frac{nq_0}{n-q_0}$, which implies
\begin{align}\label{eq2.6}
\owncount
-\frac{1}{2}-\frac{n}{2}\left(\frac{1}{q_0}-\frac{1}{\theta}\right)>-1.
\end{align}
The first inequality for $q_0$ above allows the estimation of
\begin{align*}
&\quad\ \Vert\Phi_2(u,v)(t)\Vert_{\W^{1,\theta}(\Omega)}\nonumber\\
&\leq\Vert e^{\frac{t}{\tau}(\Delta-1)}v_0\Vert_{\W^{1,\theta}(\Omega)}+\frac{1}{\tau}\intot\Vert\etsdm u(s)\Vert_{\W^{1,\theta}(\Omega)}\intd s+\frac{1}{\tau}\intot\Vert\etsdm f(s)\Vert_{\W^{1,\theta}(\Omega)}\intd s\nonumber\\
&\leq \Vert e^{\frac{t}{\tau}\Delta}v_0\Vert_{\W^{1,\theta}(\Omega)}+\frac{C_4}{\tau}\intot (t-s)^{-\frac{1}{2}-\frac{n}{2}(\frac{1}{q_0}-\frac{1}{\theta})}\Vert u(s)\Vert_{\Lo^{q_0}(\Omega)}\intd s\nonumber\\
&\hspace*{191pt}+\frac{C_5}{\tau}\intot(t-s)^{-\frac{1}{2}-\frac{n}{2}(\frac{1}{q_0}-\frac{1}{\theta})}\Vert f(s)\Vert_{\Lo^{q_0}(\Omega)}\intd s\nonumber\\
&\leq K\Vert v_0\Vert_{\W^{1,\theta}(\Omega)}+\left(\frac{C_6}{\tau}\Vert u\Vert_{\Lo^\infty(\Omega)}+\frac{C_5}{\tau}\Vert f\Vert_{\Lo^\infty([0,\infty);\,\Lo^{q_0}(\Omega))}\right)\intot(t-s)^{-\frac{1}{2}-\frac{n}{2}(\frac{1}{q_0}-\frac{1}{\theta})}\intd s\nonumber
\end{align*}
for all $t\in(0,T)$. Recalling \eqref{eq2.6} we obtain
\begin{align}\label{eq2.7}
\owncount
\Vert\Phi_2(u,v)(t)\Vert_{\W^{1,\theta}(\Omega)}\leq KR+\left(C_8(R+1)+C_7\Vert f\Vert_{\Lo^\infty([0,\infty);\,\,\Lo^{\frac{n}{2}+\delta_0}(\Omega))}\right)T^{\frac{1}{2}-\frac{n}{2}(\frac{1}{q_0}-\frac{1}{\theta})}.
\end{align}
Now if we fix $T_0\in(0,1)$ small enough and $T\in(0,T_0)$, it follows from equations \eqref{eq2.5} and \eqref{eq2.7} that $\Phi$ maps from $S$ onto itself. Furthermore using the same methods we can show that $\Phi$ acts as a contraction on $S$. To this end we pick $(u,v)\in S$ and $(\hat{u},\hat{v})\in S$ and estimate
\begin{align*}
&\quad\ \Vert \Phi_1(u,v)(t)-\Phi_1(\hat{u},\hat{v})(t)\Vert_{\Lo^\infty(\Omega)}\\
&\leq C_9\intot\Vert A^\alpha\etsd\nabla\divdot\left(u(s)\nabla v(s)-\hat{u}(s)\nabla\hat{v}(s)\right)\Vert_{\Lo^\infty(\Omega)}\intd s\\
&\leq C_{10}\intot(t-s)^{-\frac{1}{2}-\alpha-\frac{n}{2}(\frac{1}{\theta}-\frac{1}{p})}\bigg(\Vert u(s)\Vert_{\Lo^\infty(\Omega)}\Vert\nabla v(s)-\nabla \hat{v}(s)\Vert_{\Lo^\theta(\Omega)}\\
&\hspace*{221pt}
+\Vert u(s)-\hat{u}(s)\Vert_{\Lo^\infty(\Omega)}\Vert\nabla\hat{v}(s)\Vert_{\Lo^\theta(\Omega)}\bigg)\intd s\\
&\leq C_{10} T^{\frac{1}{2}-\alpha-\frac{n}{2}(\frac{1}{\theta}-\frac{1}{p})}\bigg((R+1)+(KR+1)\bigg)\bigg(\sup_{s\in[0,T]}\!\!\Vert u(s)-\hat{u}(s)\Vert_{\Lo^\infty(\Omega)}\\
&\hspace*{221pt}
+\sup_{s\in[0,T]}\!\!\Vert\nabla v(s)-\nabla\hat{v}(s) \Vert_{\Lo^\theta(\Omega)}\bigg)\\
&\leq C_{11} T^{\frac{1}{2}-\alpha-\frac{n}{2}(\frac{1}{\theta}-\frac{1}{p})}\bigg((R+1)+(KR+1)\bigg)\Vert (u,v)-(\hat{u},\hat{v})\Vert_X\ \text{ for }t\in(0,T).
\end{align*}
Analogously $\Vert \Phi_2(u,v)(t)-\Phi_2(\hat{u},\hat{v})(t)\Vert_{\W^{1,\theta}(\Omega)}\leq C_{12}T^\frac{1}{2}\Vert (u,v)-(\hat{u},\hat{v})\Vert_X\ \text{ for }t\in(0,T)$. Thus, by choosing $T\in(0,T_0)$ sufficiently small, we can ensure the contractivity of $\Phi$ on $S$ and an application of the Banach fixed-point theorem implies the existence of $(u,v)\in S$ such that $\Phi(u,v)=(u,v)$. Since $\Phi$ is in fact the variation of constants formula for the problem \eqref{KSf}, this fixed point is indeed a mild solution for \eqref{KSf}. Additionally, based on the claim preceding the arguments above, the choice of $T$ only depends on $\|u_0\|_{\Lo^{\infty}(\Omega)}$ and $\|v_0\|_{\W^{1,\theta}(\Omega)}$ and thus standard extension arguments immediately imply the existence of a maximal existence time $T_{max}$ satisfying \eqref{eq2.1}. Using standard arguments involving the semigroup estimates, as well as parabolic Schauder estimates (see \cite[Theorem IV.5.3]{lsu} for instance), it can be checked that $(u,v)$ lies in the asserted regularity class and hence solves \eqref{KSf} in the classical sense. That both $u$ and $v$ are nonnegative then follows from the comparison principle for classical sub- and supersolutions of scalar parabolic equations, by first taking $\underline{u}\equiv0$ as subsolution in the first equation and then -- since we now know both $u$ and $f$ to be nonnegative -- $\underline{v}\equiv0$ as subsolution in the second equation.

\textbf{Uniqueness: }We follow the methods of \cite{gz98} and \cite{win10ctax}. Given $T>0$ and two solutions $(u,v),(\hat{u},\hat{v})$ in $\Omega\times(0,T)$, we fix $T^*\in(0,T)$ and let $w:=u-\hat{u}$, $z:=v-\hat{v}$. Testing the difference for $(u,v)$ and $(\hat{u},\hat{v})$ of the first equation in \eqref{KSf} with $w$ we obtain for $0<t<T^*$:
\begin{align}\label{eq2.8}
\owncount
\intomega w_t w&=\intomega w\Delta w -\intomega w\nabla\divdot(u\nabla v-\hat{u}\nabla v+\hat{u}\nabla v-\hat{u}\nabla\hat{v}).\nonumber\\
\intertext{Integrating by parts und reordering of the terms we find that}\frac{1}{2}\frac{\intd}{\intd t}\intomega w^2&+\intomega\vert\nabla w\vert^2=\intomega\nabla w\cdot w\nabla v+\intomega\hat{u}\nabla z\cdot \nabla w\,\text{ for }0<t<T^*,
\intertext{where every boundary integral disappears because of the Neumann boundary conditions imposed on $u,v,\hat{u}$ and $\hat{v}$. Similar testing the second equation by $\Delta z$ results in}\label{eq2.9}\owncount
\intomega\vert\Delta z\vert^2&+\intomega\vert\nabla z\vert^2+\frac{\tau}{2}\frac{\intd}{\intd t}\intomega\vert\nabla z\vert^2=-\intomega w\Delta z\,\text{ for }0<t<T^*.
\end{align}
We continue by estimating $I_2:=\sintomega\!w\nabla v\divdot\!\nabla w,$ $I_3:=\sintomega\!\hat{u}\nabla w\divdot\!\nabla z$ and $I_4:=-\sintomega\!w\Delta z$.
Applying H\"older's inequality twice to $I_2$ we see that
\begin{align}\label{eq2.10}
\owncount
I_2&\leq\bigg(\intomega\vert\nabla w\vert^2\bigg)^\frac{1}{2}\bigg(\intomega\vert w\vert^\frac{2\theta}{\theta-2}\bigg)^\frac{\theta-2}{2\theta}\bigg(\intomega\vert\nabla v\vert^\theta\bigg)^\frac{1}{\theta}.
\intertext{By Lemma \ref{conservmass} we have $\sintomega w=0$, and thus the Gagliardo-Nierenberg inequality and the Poincaré inequality imply $
\Vert w\Vert_{\Lo^{\frac{2\theta}{\theta-2}}(\Omega)}\,\leq C_{13}\|w\|_{\W^{1,\theta}(\Omega)}^{\nicefrac{n}{\theta}}\|w\|_{\Lo^2(\Omega)}^{\nicefrac{(\theta-n)}{\theta}}\nonumber
\leq C_{14}(\,\intomega \vert\nabla w\vert^2)^\frac{n}{2\theta}(\,\intomega w^2)^\frac{\theta-n}{2\theta}$, which combined with \eqref{eq2.10} leads to}
I_2&\leq C_{14}\bigg(\intomega\vert\nabla w\vert^2\bigg)^{\frac{1}{2}+\frac{n}{2\theta}}\bigg(\intomega w^2\bigg)^{\frac{\theta-n}{2\theta}}\bigg(\intomega \vert\nabla v\vert^\theta\bigg)^\frac{1}{\theta}.\nonumber
\intertext{Finally, using the fact that $\Vert\nabla v\Vert_{\Lo^{\theta}(\Omega)}$ is bounded for $t\in(0,T^*)$ and $\theta>n\geq2$, Young's inequality from Lemma \ref{young} with $\varepsilon=\frac{1}{2C_{14}}$ and $p=\frac{1}{\frac{1}{2}+\frac{n}{2\theta}}$ yields}
I_2&\leq\frac{1}{2}\intomega\vert\nabla w\vert^2+C_{15}\intomega w^2.\nonumber
\end{align}
For $I_3$, we simply apply Young's inequality with $\varepsilon=\frac{1}{2}$ and $p=2$ to obtain
\begin{align*}
I_3\leq \frac{1}{2}\intomega\vert\nabla w\vert^2+C_{16}\intomega\vert\nabla z\vert^2\vert\hat{u}\vert^2\leq\frac{1}{2}\intomega\vert\nabla w\vert^2+C_{17}\intomega\vert\nabla z\vert^2
\end{align*}
by the boundedness of $\hat{u}$ in $\Omega\times(0,T^*)$. In a similar fashion
\begin{align*}
I_4\leq\intomega\vert\Delta z\vert^2+\frac{1}{4}\intomega w^2.
\end{align*}
Adding up \eqref{eq2.8} and \eqref{eq2.9} and using the above estimates for $I_2,I_3$ and $I_4$, we have
\begin{align*}
&\frac{1}{2}\frac{\intd}{\intd t}\intomega w^2+\intomega\vert\nabla w\vert^2+\intomega\vert\Delta z\vert^2+\intomega\vert\nabla z\vert^2+\frac{\tau}{2}\frac{\intd}{\intd t}\intomega\vert\nabla z\vert^2=I_2+I_3+I_4\\
\leq\ &\frac{1}{2}\intomega\vert\nabla w\vert^2+C_{18}\intomega w^2+\frac{1}{2}\intomega\vert\nabla w\vert^2+C_{18}\intomega\vert\nabla z\vert^2+\intomega\vert\Delta z\vert^2+\frac{1}{4}\intomega w^2
\end{align*}
for $t\in(0,T^*)$. Therefore the differential inequality
\begin{align*}
\frac{\intd}{\intd t}\bigg(\intomega w^2+\tau\intomega\vert\nabla z\vert^2\bigg)\leq C_{19}\bigg(\intomega w^2+\tau\intomega\vert\nabla z\vert^2\bigg)
\end{align*}
holds for all $t\in(0,T^*)$, which implies $w\equiv0$ and $z\equiv0$ in $\Omega\times(0,T^*)$. Since $T^*<T_{max}$ was arbitrary, the claim follows.
\end{proof}

We should mention that the proof above shows that assuming \eqref{thetareq} always to be true in the following sections is no restriction. This will be required to appropriately control the powers appearing in the integrals required for several estimates.
%

\setcounter{gleichung}{0} 

\section{Boundedness of solutions in the two dimensional case}\label{cglob2}
%
After establishing time local existence our next objective is to prove, under certain conditions, the existence of time-independent bounds for the norms $\|u\|_{\Lo^\infty(\Omega)}$ and $\|v\|_{\W^{1,\theta}(\Omega)}$ with $\theta$ satisfying $n<\theta<\frac{n^2+2n\delta_0}{n-2\delta_0}$. These time-independent bounds then also imply time global existence by means of \eqref{eq2.1}.

As mentioned in the introduction we will first limit ourselves to the two dimensional case. To this end we assume throughout this section $n=2$, that $\theta$ satisfies the inequalities $2<\theta<\frac{2+2\delta_0}{1-\delta_0}$, and that all other necessary conditions of Theorem \ref{localexist} are fulfilled. Accordingly we denote by $(u,v)$ the classical solution to \eqref{KSf} in dimension $n=2$, given by Theorem \ref{localexist}.
%

\subsection{Essential estimates}\label{cglob2s1}
Let us first gather some essential results needed for the proofs later on.
To obtain our main result for this section, we will roughly follow the methods of \cite{NSY97}. One essential ingredient to this method will be Moser's sharp form of an inequality by Trudinger, which can be found in \cite{trud67} and \cite{mos71} respectively. This result was extended by Nagai, Senba and Yoshida in \cite{NSY97} for the two-dimensional radial symmetric case, and by Chang and Yang in \cite{CY88} without the restriction to radially symmetric functions. Since the methods surrounding these inequalities in our proofs are identical in both cases, we only focus on general domains and thus only give the result by \cite{CY88} adjusted to our notation.
\begin{theo}[A Trudinger-Moser type inequality]\label{trudingermoser}
Let $\Omega$ be a smoothly bounded domain in $\mathbb{R}^2$. Then there exists a constant $C>0$ depending on $\Omega$ such that
\begin{align*}
\intomega e^{\vert v\vert}\intd x\leq C e^{\left(\frac{1}{8\pi}\|\nabla v\|^2_{\Lo^2(\Omega)}+\frac{1}{\vert\Omega\vert}\|v\|_{\Lo^1(\Omega)}\right)}
\end{align*}
holds for all $v\in\W^{1,2}(\Omega)$.
\end{theo}
Furthermore we will make use of the following inequality, which is a modification of an inequality found in \cite{bil94}.
\begin{lem}\label{f-biler-ineq}
Let $\Omega$ be a smoothly bounded domain in $\mathbb{R}^2$ and $2\leq p<\infty$. Then for any $\varepsilon>0$ there exists a constant $C>0$ depending only on $p$ and $\varepsilon$, such that
\begin{align*}
\|v\|_{\Lo^p(\Omega)}\leq\varepsilon\|\nabla v\|_{\Lo^2(\Omega)}^{1-\nicefrac{1}{p}}\|v\log|v|\|_{\Lo^1(\Omega)}^{\nicefrac{1}{p}}+C\left(\|v\log|v|\|_{\Lo^1(\Omega)}+\|v\|^{\nicefrac{1}{p}}_{\Lo^1(\Omega)}+\|v\|_{\Lo^1(\Omega)}\right)
\end{align*}
holds for every $v\in\W^{1,2}(\Omega)$.
\end{lem}
\begin{bew} The proof closely follows arguments shown in \cite{bil94} or \cite[Lemma A.5]{TW14} and will only be given for  convenience of the reader. Consider $N>1$, to be specified below, and the piecewise defined function $\chi:\mathbb{R}\rightarrow\mathbb{R}$ given by
\begin{align*}
\chi(s)=\begin{cases}
0\quad&,\text{ if }\vert s\vert\leq N\\
2\vert s\vert-2N&,\text{ if }N< \vert s\vert< 2N\\
\vert s\vert&,\text{ if }2N\leq\vert s\vert.
\end{cases}
\end{align*}
Writing $Y_1:=\lbrace x\in\Omega\ \vert\ \vert v\vert\leq N\rbrace$ and $Y_2:=\lbrace x\in\Omega\ \vert\  N<\vert v\vert< 2N\rbrace$, we can estimate
\begin{align}
\left\| \vert v\vert -\chi(v)\right\|_{\Lo^p(\Omega)}^p&=\int\limits_{Y_1}\vert v\vert^p\intd x+\int\limits_{Y_2}\left\vert 2N-\vert v\vert\right\vert^p\intd x\nonumber\leq N^{p-1}\int\limits_{Y_1}\vert v\vert\intd x+(2N)^{p-1}\!\int\limits_{Y_2}\vert 2N-\vert v\vert\vert\intd x.\nonumber
\intertext{Since $N<\vert v\vert<2N$ implies $0<2N-\vert v\vert<\vert v\vert$, the second integral can be further estimated to obtain}\label{fbilerproof-eq1}\owncount
\left\| \vert v\vert -\chi(v)\right\|_{\Lo^p(\Omega)}^p\,&\leq N^{\,p-1}\!\!\int\limits_{Y_1}\vert v\vert\intd x+(2N)^{p-1}\!\int\limits_{Y_2}\vert v\vert\intd x\leq\left(2N\right)^{p-1}\!\!\!\int\limits_{Y_1\cup\,Y_2}\!\!\!\vert v\vert\intd x\leq\left(2N\right)^{p-1}\|v\|_{\Lo^1(\Omega)}.
\end{align}
Additionally, using that by the definition of $\chi$ there holds $\chi(v)\leq\vert v\vert$ for all $x\in\Omega$ and defining the set $Y_3:=\lbrace x\in\Omega\ \vert\ \vert v\vert\geq N \rbrace$, we estimate
\begin{align}
\owncount\label{fbilerproof-eq2}
\|\chi(v)\|^{\nicefrac{1}{p}}_{\Lo^1(\Omega)}&=\left(\ \int\limits_{Y_3}\chi(v)\intd x\right)^{\frac{1}{p}}\!\leq\left(\ \int\limits_{Y_3} \vert v\vert\frac{\log\vert v\vert}{\log(N)}\intd x\right)^{\frac{1}{p}}\leq\log(N)^{-\frac{1}{p}}\|v\log\vert v\vert\|_{\Lo^1(\Omega)}^{\nicefrac{1}{p}},
\end{align}
as well as
\begin{align}\label{fbilerproof-eq3}
\owncount
\|\chi(v)\|_{\W^{1,2}(\Omega)}^{1-\nicefrac{1}{p}}\,&=\left(\|\nabla(\chi(v))\|_{\Lo^2(\Omega)}^2+\|\chi(v)\|_{\Lo^2(\Omega)}^2\right)^{\frac{p-1}{2p}}\leq\left(\|\chi'(v)\nabla v\|_{\Lo^2(\Omega)}^2+\|v\|_{\Lo^2(\Omega)}^2\right)^{\frac{p-1}{2p}}\nonumber\\
&\leq\left(4\|\nabla v\|_{\Lo^2(\Omega)}^2+\|v\|_{\Lo^2(\Omega)}^2\right)^{\frac{p-1}{2p}}\leq 2^{\frac{p-1}{p}}\|v\|_{\W^{1,2}(\Omega)}^{1-\nicefrac{1}{p}}.
\end{align}
Using the Gagliardo-Nierenberg inequality in combination with \eqref{fbilerproof-eq1}, \eqref{fbilerproof-eq2} and \eqref{fbilerproof-eq3} we can estimate  $\|v\|_{\Lo^p(\Omega)}$ in the following way:
\begin{align*}
\|v\|_{\Lo^p(\Omega)}\,&\leq\|\vert v\vert-\chi(v)\|_{\Lo^p(\Omega)}+\|\chi(v)\|_{\Lo^p(\Omega)}\\
&\leq (2N)^\frac{p-1}{p}\|v\|_{\Lo^1(\Omega)}^{\nicefrac{1}{p}}+C_{1}\|\chi(v)\|_{\W^{1,2}(\Omega)}^{1-\nicefrac{1}{p}}\|\chi(v)\|_{\Lo^1(\Omega)}^{\nicefrac{1}{p}}\\
&\leq(2N)^\frac{p-1}{p}\|v\|_{\Lo^1(\Omega)}^{\nicefrac{1}{p}}+\frac{2^{\frac{p-1}{p}}C_{1}}{\log(N)^\frac{1}{p}}\|v\|_{\W^{1,2}(\Omega)}^{1-\nicefrac{1}{p}}\|v\log\vert v\vert\|_{\Lo^1(\Omega)}^{\nicefrac{1}{p}}.
\end{align*}
Now we apply the Poincaré inequality to obtain
\begin{align*}
\|v\|_{\Lo^p(\Omega)}&\leq(2N)^\frac{p-1}{p}\|v\|_{\Lo^1(\Omega)}^{\nicefrac{1}{p}}+\frac{C_{2}}{\log(N)^\frac{1}{p}}\left(\|\nabla v\|_{\Lo^2(\Omega)}^{1-\nicefrac{1}{p}}+\|v\|_{\Lo^1(\Omega)}^{1-\nicefrac{1}{p}}\right)\|v\log\vert v\vert\|_{\Lo^1(\Omega)}^{\nicefrac{1}{p}}.
\intertext{Finally, using Young's inequality and taking $N=exp\left(\left(\frac{C_2}{\varepsilon}\right)^p\right)$ results in}
\|v\|_{\Lo^p(\Omega)}&\leq\varepsilon\|\nabla v\|_{\Lo^2(\Omega)}^{1-\nicefrac{1}{p}} \|v\log |v|\|_{\Lo^1(\Omega)}^{\nicefrac{1}{p}}+C(p,\varepsilon)\left(\|v\log |v|\|_{\Lo^1(\Omega)}+\|v\|_{\Lo^1(\Omega)}^{\nicefrac{1}{p}}+\|v\|_{\Lo^1(\Omega)}\right)
\end{align*}
and thus completing the proof.
\end{bew}
The last preparatory result we take advantage of is in fact a simple consequence of Young's inequality.
\begin{lem}\label{Malpha-ineqlem}
Assume $C_1>0,C_2>0$ and $\beta\in(0,1)$ to be constant. Then there exists a constant $M_0>0$, such that every $M\geq0$ satisfying
\begin{align}\label{malphaeq0}
\owncount
M\leq C_1+C_2 M^\beta
\end{align}
fulfills the inequality $M\leq M_0$.
\end{lem}
\begin{bew}
We apply Young's inequality with $\varepsilon:=\frac{1}{2C_2},a:=M^\beta, b:=1,p:=\frac{1}{\beta}$ and $q:=\frac{p}{p-1}$ in \eqref{malphaeq0} to obtain
\begin{align*}
2M\leq&\ 2\left(C_1+C_2M^\beta\right)\leq\ 2C_1+2C_2\left(\varepsilon a^p+C(\varepsilon,p,q)b^{\frac{p}{p-1}}\right).
\end{align*}
With $C(p,q):=C(\varepsilon,p,q)\cdot\varepsilon^{\frac{q}{p}}=p^{-\frac{q}{p}}q^{-1}=\frac{p-1}{p^{\frac{p}{p-1}}}$ this yields
\begin{align}\label{malphaeq1}
\owncount
2M\leq 2C_1+2C_2\left(\frac{1}{2C_2}M+C(p,q)\varepsilon^{-\frac{1}{p-1}}\right)=2C_1+M+2C_2 C(p,q)(2C_2)^{\frac{\beta}{1-\beta}}.
\end{align}
Thus reordering the terms in \eqref{malphaeq1} results in the asserted inequality for $M$.
\end{bew}

Before we can start our energy-based approach, we require some lemmata regarding the $\Lo^p(\Omega)$-norms of the chemical concentration function $v$ and the respective gradient function $\nabla v$. The lemma below can not only be proven for $n=2$, but also in general dimensions. The results follows almost immediately from integrating the second equation of \eqref{KSf} and thus will be omitted here.
\begin{lem}\label{massforv}
For all $t\in[0,T_{max})$ it holds that:
\begin{align*}
\|v(\cdot,t)\|_{\Lo^1(\Omega)}&=e^{-\frac{t}{\tau}}\|v_0\|_{\Lo^1(\Omega)}+\|u_0\|_{\Lo^1(\Omega)}(1-e^{-\frac{t}{\tau}})+\frac{1}{\tau}\intot\|f(\cdot,s)\|_{\Lo^1(\Omega)}e^{\frac{s-t}{\tau}}\intd s.
\intertext{In particular there is some $C>0$ such that}
\|v(\cdot,t)\|_{\Lo^1(\Omega)}&\leq\,\|v_0\|_{\Lo^{1}(\Omega)}+\|u_0\|_{\Lo^{1}(\Omega)}+C\|f\|_{\Lo^\infty([0,\infty);\,\Lo^{1+\delta_0}(\Omega))}
\end{align*}
holds for all $t\in[0,T_{max})$.
\end{lem}

The first step in a quite long estimation process will be the extraction of $L^p(\Omega)$ bounds for $v$ and $\nabla v$ for some $p\in[1,\infty]$, by making use of appropriate semigroup estimates.
\begin{prop}\label{v estimate}
Let $q\in[1,2)$ and $p\in[1,\infty)$. Then there exists a time-independent constant $C>0$ such that
\begin{align*}
\|\nabla v(\cdot,t)\|_{\Lo^q(\Omega)}&\leq C\ \text{ and }\ \|v(\cdot,t)\|_{\Lo^p(\Omega)}\leq C,
\end{align*}
hold for all $0<t<T_{max}$.
\end{prop}
\begin{bew}
Since $q<2<\theta$, using the semigroup estimates again, we find time independent constants $C_i>0$ such that
\pagebreak
\begin{align*}
\|\nabla v(\cdot,t)\|_{\Lo^q(\Omega)}&\leq C_1\|\nabla v_0\|_{\Lo^\theta(\Omega)}+C_2\!\intot\left(1+(t-s)^{-\frac{1}{2}-(1-\frac{1}{q})}\right)e^{-(\lambda_1+1)(t-s)}\|u(\cdot,s)\|_{\Lo^1(\Omega)}\intd s\\
&\hspace*{70pt}+C_3\!\intot\left(1+(t-s)^{-\frac{1}{2}-(1-\frac{1}{q})}\right)e^{-(\lambda_1+1)(t-s)}\|f(\cdot,s)\|_{\Lo^1(\Omega)}\intd s
\intertext{for $t\in[0,T_{max})$, where $\|f(\cdot,t)\|_{\Lo^1(\Omega)}\leq C_4\|f\|_{\Lo^\infty([0,\infty);\,\Lo^{1+\delta_0}(\Omega))}$ and $\|u(\cdot,t)\|_{\Lo^1(\Omega)}=\|u_0\|_{\Lo^1(\Omega)}$ by Lemma \ref{conservmass}. Thus, using the fact that $-\frac{1}{2}-(1-\frac{1}{q})>-1$ for $q<2$, we obtain}
\|\nabla v(\cdot,t)\|_{\Lo^q(\Omega)}&\leq C_5+C_6\intot\left(1+(t-s)^{-\frac{1}{2}-(1-\frac{1}{q})}\right)e^{-(\lambda_1+1)(t-s)}\intd s\\
& \leq C_5+C_6\int\limits_0^\infty\!\!\left(1+\sigma^{-\frac{1}{2}-(1-\frac{1}{q})}\right)e^{-(\lambda_1+1)\sigma}\intd\sigma=C_5+C_6\!\left(\frac{1}{\lambda_1+1}+\frac{\Gamma(\frac{1}{q}-\frac{1}{2})}{(\lambda_1+1)^{\frac{1}{q}-\frac{1}{2}}}\right)
\end{align*}
for $t\in[0,T_{max})$, with the Gamma-function $\Gamma(x)=\int_0^\infty t^{x-1}e^{-t}\intd t$, proving the first asserted inequality. To verify the second inequality we apply the Poincaré inequality in combination with the Sobolev embedding theorem to obtain $$\|v\|_{\Lo^{\frac{2q}{2-q}}(\Omega)}\leq C_7\|v\|_{\W^{1,q}(\Omega)}\leq C_8\|\nabla v\|_{\Lo^q(\Omega)}+C_8\|v\|_{\Lo^1(\Omega)}$$ for every $q<2$ and $t\in[0,T_{max})$, which by Lemma \ref{massforv} completes the proof.
\end{bew}

Thanks to H\"older's inequality and the nonnegativity of $f$ and $v$, an immediate consequence is a bound for $\sintomega fv\intd x$.
\begin{cor}\label{fvbound}
There exists a constant $C>0$, such that $\sintomega fv\intd x\leq C$ holds for all $t\in[0,T_{max})$.
\end{cor}
\subsection{Boundedness in two dimensions~- an energy-based approach}\label{cglob2s2}
In addition to all the assumptions mentioned at the beginning of Section \ref{cglob2}, we now also assume $f$ to be constant-in-time.

The next step in our estimation process involves defining a suitable energy functional of the system \eqref{KSf}. For the system without external signal production the corresponding functional can be defined as $W'(u,v):=\sintomega\left[u \log u -uv+\frac{1}{2}(\vert\nabla v\vert^2+v^2)\right]\intd x$ (cf. \cite{NSY97}). Using the derivation method of \cite{NSY97}, it follows that for nonnegative and constant-in-time signal production $f$ this functional naturally extends to $$W(u,v):=\intomega\left[u \log u -uv+\frac{1}{2}(\vert\nabla v\vert^2+v^2)-fv\right]\intd x.$$ In particular we obtain the following
\begin{lem}\label{energylem}
Let $(u,v)$ be the classical solution of \eqref{KSf} in $\Omega\times(0,T_{max})$ and
\begin{align*}
W(t):=W\left(u(\cdot,t),v(\cdot,t)\right)=\intomega\left[u\log u-uv+\frac{1}{2}(\vert\nabla v\vert^2+v^2)-fv\right]\intd x.
\end{align*}
Then
\begin{align*}
\frac{\intd}{\intd t}W(t)+\intomega u\vert\nabla(\log u-v)\vert^2\intd x+\tau\intomega v_t^2\intd x=0.
\end{align*}
holds for all $t\in(0,T_{max})$.
\end{lem}
\begin{bew}
Writing the first equation of \eqref{KSf} as $u_t=\nabla\divdot(\nabla u-u\nabla v)$ and testing this with $(\log u-v)$ yields
\begin{align*}
\intomega u_t(\log u-v)\intd x\,&=\intomega(\log u -v)\nabla\divdot(\nabla u-u\nabla v)\intd x=-\intomega(\nabla u-u\nabla v)\divdot\nabla(\log u-v)\intd x\nonumber
\end{align*}
for all $t\in(0,T_{max})$. Since $u\nabla(\log u-v)=\nabla u-u\nabla v$, this results in
\begin{align}\label{energylemeq1}
\owncount
\intomega u_t(\log u-v)\intd x\ &=-\intomega u\vert\nabla(\log u-v)\vert^2\intd x\ \text{for all }t\in(0,T_{max}).
\end{align}
Furthermore rewriting the left side of the equation as
\begin{align*}
\intomega u_t(\log u-v)\intd x\,=\intomega\left[ (u\log u)_t-u_t-u_t v\right]\intd x
=\frac{\intd}{\intd t}\intomega u\log u\intd x-\frac{\intd}{\intd t}\intomega u\intd x-\intomega u_t v\intd x
\end{align*}
for all $t\in(0,T_{max})$, and using Lemma \ref{conservmass} to deduce $\frac{\intd}{\intd t}\sintomega u\intd x=0$, we obtain
\begin{align}\label{energylemeq2}
\owncount
\intomega u_t(\log u-v)\ &=\frac{\intd}{\intd t}\intomega u\log u\intd x-\frac{\intd}{\intd t}\intomega uv\intd x+\intomega v_t u\intd x\ \text{for all }t\in(0,T_{max}).
\end{align}
Utilizing the second equation of \eqref{KSf} to express $u$, it holds for $t\in(0,T_{max})$ that
\begin{align}\label{energylemeq3}
\owncount
\intomega v_t u\intd x\,&=\tau\intomega v_t^2\intd x+\intomega v v_t\intd x-\intomega \Delta v v_t\intd x-\intomega f v_t\intd x\nonumber\\
&=\tau\intomega v_t^2\intd x+\frac{1}{2}\frac{\intd}{\intd t}\intomega v^2\intd x+\frac{1}{2}\frac{\intd}{\intd t}\intomega\vert\nabla v\vert^2\intd x-\frac{\intd}{\intd t}\intomega fv\intd x,
\end{align}
where we used $f_t\equiv0$ and $-\sintomega\Delta v v_t\intd x=\sintomega\,(\frac{\intd}{\intd t}\nabla v)\divdot\!\nabla v\intd x$. Combining \eqref{energylemeq1},\eqref{energylemeq2} and \eqref{energylemeq3} we have	
\begin{align*}
-\intomega u\vert\nabla(\log u -v)\vert^2\intd x=\frac{\intd}{\intd t}\intomega\left[u\log u-uv+\frac{1}{2}(v^2+\vert\nabla v\vert^2)-fv\right]\intd x+\tau\intomega v_t^2\intd x
\end{align*}
for all $t\in(0,T_{max})$, and thus by definition of $W(t)$ the proof is complete.
\end{bew}
An immediate consequence of the lemma above is the following
\begin{rem}\label{wboundrem}
Let $W$ be as in Lemma \ref{energylem}. Then $\frac{\intd}{\intd t}W(t)\leq 0$ for all $t\in(0,T_{max})$. In particular
\begin{align*}
W(t)\leq W(0)\leq\intomega u_0\log u_0+\frac{1}{2}(\vert\nabla v_0\vert^2+v_0^2)\intd x\leq C\ \text{for all }t\in[0,T_{max})
\end{align*}
for some $C>0$, by the conditions imposed on the initial values.
\end{rem}

With $W$ defined in Lemma \ref{energylem}, representing the time evolution of some type of energy, and the inequalities above we are now able to prove time-independent bounds for $\sintomega uv\intd x$ and $\vert W\vert$, which are the first main step on the way to the boundedness result.
\begin{prop}\label{criticalmassprop}
Let $W$ be defined as in Lemma \ref{energylem}. If $\sintomega u_0\intd x<4\pi$, then there exists a positive time-independent constant $C$ such that
\begin{align*}
\qquad\quad\intomega uv\intd x\leq C\ \text{ and }\ \vert W(t)\vert\leq C\quad\text{holds for all }t\in[0,T_{max}).
\end{align*}
\end{prop}
\begin{bew}
We make use of the method shown in \cite[Lemma 3.4]{NSY97}. To this end we split $W(t)$ for $\delta>0$ according to
\begin{align*}
W(t)\,&=\intomega\left[u \log u - (1+\delta)uv\right]\intd x+\intomega[\delta uv+\frac{1}{2}(\vert\nabla v\vert^2+v^2)-fv]\intd x\ \text{for all }t\in[0,T_{max}).
\end{align*}
By straightforward rearrangement we can further rewrite this equality as
\begin{align}\label{criticalmasspropeq1}
\owncount
W(t)=-\intomega u\log\left(\frac{e^{(1+\delta)v}}{u}\right) \intd x+\intomega[\delta uv+\frac{1}{2}(\vert\nabla v\vert^2+v^2)-fv]\intd x\ \text{for all }t\in[0,T_{max}).
\end{align}
Set $m_u:=\sintomega u_0\intd x$, then by Lemma \ref{conservmass} $\sintomega\frac{u}{m_u}\intd x=1$ for all $t\in[0,T_{max})$. Using the convexity of the function $-\log x$, we can apply Jensen's inequality to obtain
\begin{align}\label{criticalmasspropeq2}
\owncount
&\ -\log\left(\frac{1}{m_u}\intomega e^{(1+\delta)v}\intd x\right)\leq-\frac{1}{m_u}\intomega u\log\left(\frac{e^{(1+\delta)v}}{u}\right)\intd x\ \text{for all }t\in[0,T_{max}).
\end{align}
Invoking the Trudinger-Moser inequality from Theorem \ref{trudingermoser}, we can estimate
\begin{align}\label{criticalmasspropeq3}
\owncount
\intomega e^{(1+\delta)v}\intd x&\leq Ce^{\left(\frac{(1+\delta)^2}{8\pi}\|\nabla v\|^2_{\Lo^2(\Omega)}+\frac{(1+\delta)}{\vert\Omega\vert}\|v\|_{\Lo^1(\Omega)}\right)}\nonumber
\intertext{for all $t\in[0,T_{max})$ with some $C>0$, which by the monotonicity of the logarithm implies}
\log\left(\frac{1}{m_u}\intomega e^{(1+\delta)v}\intd x\right)&\leq\log\left(\frac{C}{m_u}\right)+\frac{(1+\delta)^2}{8\pi}\|\nabla v\|_{\Lo^2(\Omega)}^2+\frac{(1+\delta)}{\vert\Omega\vert}\|v\|_{\Lo^1(\Omega)}
\end{align}
for all $t\in[0,T_{max})$. Consequently, combining \eqref{criticalmasspropeq1} -- \eqref{criticalmasspropeq3}, we see that
\begin{align}\label{criticalmasspropeq4}
\owncount
W(t)&\geq -m_u\log\left(\frac{1}{m_u}\intomega e^{(1+\delta)v}\intd x\right)+\intomega\left(\delta uv+\frac{1}{2}(\vert\nabla v\vert^2+v^2)-fv\right)\intd x\nonumber\\
&\geq-m_u\left(\log\left(\frac{C}{m_u}\right)+\frac{(1+\delta)^2}{8\pi}\|\nabla v\|^2_{\Lo^2(\Omega)}+\frac{(1+\delta)}{\vert\Omega\vert}\|v\|_{\Lo^1(\Omega)}\right)\\\nonumber&\hspace*{90pt}+\intomega\left(\delta uv+\frac{1}{2}(\vert\nabla v\vert^2+v^2)-fv\right)\intd x\ \text{for all }t\in[0,T_{max}).
\end{align}
Reordering the terms in \eqref{criticalmasspropeq4}, we have
\begin{align}\label{criticalmasspropeq5}
\owncount
\mbox{ }&\left(\frac{1}{2}-\frac{m_u(1+\delta)^2}{8\pi}\right)\|\nabla v\|_{\Lo^2(\Omega)}^2+\delta\| uv\|_{\Lo^1(\Omega)}+\frac{1}{2}\|v\|_{\Lo^2(\Omega)}^2\nonumber\\ 
\leq&\ W(t)+m_u\left[\log\left(\frac{C}{m_u}\right)\!+\frac{(1+\delta)}{\vert\Omega\vert}\|v\|_{\Lo^1(\Omega)}\right]\!+\|fv\|_{\Lo^1(\Omega)}\quad \text{for all }t\in[0,T_{max}).
\end{align}
The bound imposed on $m_u$ implies $\frac{m_u}{8\pi}<\frac{1}{2}$, thus by choosing $\delta>0$ small enough, we have $\frac{1}{2}-\frac{m_u(1+\delta)^2}{8\pi}>0$. Since the left side of \eqref{criticalmasspropeq5} is nonnegative, all that is left to prove is the boundedness of the right side of \eqref{criticalmasspropeq5}. The bound in question follows from Lemma \ref{massforv}, Corollary \ref{fvbound} and Remark \ref{wboundrem}.
\end{bew}

From now on we assume $\sintomega u_0\intd x<4\pi$ for the rest of this section. We can quite easily extract two further bounds from the proposition above and the definition of $W(t)$ (see also \cite[Remark 3.10]{NSY97}):
\begin{cor}\label{vtrem}
There exist constants $C_1,C_2>0$, not depending on $t$, such that
\begin{align*}
\intot\|v_t(\cdot,s)\|_{\Lo^2(\Omega)}^2\intd s\leq C_1\text{ and }\|u\log u\|_{\Lo^1(\Omega)}\leq C_2
\end{align*}
holds for all $0\leq t<T_{max}$. 
\end{cor}

Recalling the alternative \eqref{eq2.1} regarding the maximal existence time from Theorem \ref{localexist}, we want to show the boundedness of the involved quantities to deduce $T_{max}=\infty$. Since the possibility of immediately estimating the required norms without further preparation seems farfetched, we continue with the more reasonable estimation of $\|u(\cdot,t)\|_{\Lo^2(\Omega)}$.

\begin{prop}\label{u2-bound}
There exists a constant $C>0$ such that $\|u(\cdot,t)\|_{\Lo^2(\Omega)}\leq C$ holds for all $t\in(0,T_{max})$.
\end{prop}
\begin{bew}
The proof follows the argumentation shown in \cite[Lemma 3.6]{NSY97} very closely. We therefore only briefly talk about additional estimation needed and the main steps necessary. For more in-depth reasoning we refer to \cite{NSY97}. Again we denote by $C_i$ positive constants independent of $t$. 

Using H\"older's inequality we have $\sintomega u^2f\intd x\leq C_1\|u\|_{\Lo^r(\Omega)}^2,$ with $r:=\frac{2(1+\delta_0)}{\delta_0}<\infty$. Additional application of the Gagliardo-Nierenberg inequality in conjunction with the Poincaré inequality and Lemma \ref{conservmass} yields $\sintomega u^2f\intd x \leq C_2\|\nabla u\|_{\Lo^2(\Omega)}^{2a}+C_3,$ where $a:=1-\frac{1}{r}<1$. Now Young's inequality with $p=\frac{1}{a}$ and small $\varepsilon>0$, to be specified below, allows us to further estimate
\begin{align}
\owncount\label{u2-bound-eq5}
\intomega u^2f\intd x&\leq \frac{\varepsilon}{2}\|\nabla u\|_{\Lo^2(\Omega)}^2+C_{4}\ \text{ for all }t\in(0,T_{max}).
\end{align}
Testing the first equation from \eqref{KSf} with $u$ and inserting the second equation of \eqref{KSf} implies
\begin{align}\label{u2-bound-eq2}
\owncount
\frac{1}{2}\frac{\intd}{\intd t}\intomega u^2\intd x+\intomega\vert\nabla u\vert^2\intd x
\leq\frac{\tau}{2}\intomega\vert u^2 v_t\vert\intd x+\frac{1}{2}\intomega u^3\intd x+\frac{1}{2}\intomega u^2f\intd x\ \text{ for all }t\in(0,T_{max}).
\end{align}
Making use of \eqref{u2-bound-eq5} and arguments involving Lemma \ref{f-biler-ineq}, the H\"older and Gagliardo-Nierenberg inequalities, as well as the Jensen and Young inequalities (see \cite[Lemma 3.6]{NSY97}) we see that
\begin{align}
\owncount\label{u2-bound-eq6}
\ &\ \frac{\intd}{\intd t}\|u\|_{\Lo^2(\Omega)}^2+\left(2-\varepsilon-\varepsilon^3\|u\log u\|_{\Lo^1(\Omega)}\right)\|\nabla u\|^2_{\Lo^2(\Omega)}\nonumber\\
\leq&\ C_{5}\left((\|v_t\|_{\Lo^2(\Omega)}^2+\|v_t\|_{\Lo^2(\Omega)})\|u\|_{\Lo^2(\Omega)}^2+\|u\log u\|_{\Lo^1(\Omega)}^3+\|u\|_{\Lo^1(\Omega)}+\|u\|_{\Lo^1(\Omega)}^3+1\right)
\end{align}
for all $t\in(0,T_{max})$. Recalling that $\|u\log u\|_{\Lo^1(\Omega)}$ is bounded -- see Corollary \ref{vtrem} -- it is possible to choose $\varepsilon>0$ such that $2-\varepsilon-\varepsilon^3\|u\log u\|_{\Lo^1(\Omega)}\geq1$. Using the Gagliardo-Nierenberg inequality to control $\|\nabla u\|_{\Lo^2(\Omega)}^2$ from below, we may express \eqref{u2-bound-eq6} as the ordinary differential inequality
\begin{align*}
g'(t)+g(t)\leq C_{6}\left(h(t)+\sqrt{h(t)}\right)g(t)+K\ \text{ for all }t\in(0,T_{max}),
\end{align*}
where we put $K:=C_{7}\sup_{t\in[0,T_{max})}\!\left(1+\!\|u\log u\|_{\Lo^1(\Omega)}^3\!+\|u\|_{\Lo^1(\Omega)}^2+\!\|u\|_{\Lo^1(\Omega)}+\!\|u\|_{\Lo^1(\Omega)}^3\right)$, $h(t):=\|v_t(\cdot,t)\|_{\Lo^2(\Omega)}^2$ and $g(t):=\|u(\cdot,t)\|_{\Lo^2(\Omega)}^2$.
Since $(C_{6}\sqrt{h(t)}-1)^2\geq0$, it is $C_{6}\sqrt{h(t)}\leq\frac{1}{2}+\frac{1}{2}C_{6}^2h(t)$ and thus setting $C_8:=C_{6}+\frac{1}{2}C_{6}^2$ we obtain
\begin{align*}
g'(t)+\left(\frac{1}{2}-C_8 h(t)\right)g(t)\leq K\ \text{ for all }t\in(0,T_{max}).
\end{align*}
Consequently, defining $\varphi(t):=\intot\left(\frac{1}{2}-ch(s)\right)\intd s$ -- which in view of Corollary \ref{vtrem} is bounded from below by $\frac{1}{2}t-C_{9}$ -- this implies
\begin{align*}
\|u(\cdot,t)\|_{\Lo^2(\Omega)}^2\leq C_{10}\left(\|u_0\|_{\Lo^2(\Omega)}^2 e^{-\frac{t}{2}}+K\right)\leq C_{11}\ \text{ for all }t\in(0,T_{max}),
\end{align*}
completing the proof.
\end{bew}

Using the bound above we can extend the regularity result for $\nabla v$ to the space $\Lo^\theta(\Omega)$ with $\theta$ satisfying $2<\theta<\frac{2+2\delta_0}{1-\delta_0}$. In fact we have

\begin{prop}\label{propnablavtheta}
Let $\theta$ satisfy the inequalities $2<\theta<\frac{2+2\delta_0}{1-\delta_0}$. Then there exists a constant $C_1>0$ such that $\|\nabla v(\cdot,t)\|_{\Lo^\theta(\Omega)}\leq C_1$ holds for all $t\in(0,T_{max})$. In particular, $\|v(\cdot,t)\|_{\W^{1,\theta}(\Omega)}\leq C_2$ for all $t\in(0,T_{max})$ and some suitable constant $C_2>0$.  
\end{prop}
\begin{bew}
Apply the semigroup estimates for $\|\nabla v\|_{\Lo^\theta(\Omega)}$ and make use of the bounds shown in Proposition \ref{u2-bound}, Lemma \ref{conservmass} and Lemma \ref{massforv} as well as the Poincaré inequality.
\end{bew}
Additionally, the proposition above was the last missing piece in order to establish 

\begin{prop}\label{propuinf}
There exists a constant $C>0$ such that $\|u(\cdot,t)\|_{\Lo^\infty(\Omega)}\leq C$ for all $t\in(0,T_{max})$. 
\end{prop}
\begin{bew} A similar, and for even more general Keller-Segel systems applicable, argumentation to the one given below, can be found in \cite[Lemma 3.2]{BBWT15}. Nevertheless we give a full proof here for the sake of completion.\\
Let $A:=-\Delta+1$. Next we choose $p\in(2,\min\lbrace\theta,4\rbrace)$ and $q>\frac{2p}{p-2}$. These choices imply $q>p$ and that the interval $I:=(\frac{2}{q},\frac{1}{2}-\frac{1}{p}+\frac{1}{q})$ is not empty, which in turn allows us to choose $\alpha\in I$ satisfying $q\alpha>2$ as well as $0<\alpha<1$. Thus we can find a constant $C>0$ such that the fractional power $A^\alpha$ fulfills $\|w\|_{\Lo^\infty(\Omega)}\leq C\|A^\alpha w\|_{\Lo^q(\Omega)}$ and $\|A^\alpha e^{\sigma\Delta}w\|_{\Lo^q(\Omega)}\leq C\sigma^{-\alpha}\|w\|_{\Lo^q(\Omega)}$ for all $w\in C^{\infty}_0(\Omega)$ and $\sigma>0$ (cf. \cite[1.6.1 and 1.4.3]{hen81}). Furthermore these choices imply $-\frac{1}{2}-\alpha-(\frac{1}{p}-\frac{1}{q})>-1$. Combining these statements with semigroup estimates and the maximum principle leads to
\begin{align}\label{propuinf-eq1}
\owncount\mbox{ }&\ \|u(\cdot,t)\|_{\Lo^\infty(\Omega)}\\
\leq&\ \|u_0\|_{\Lo^\infty(\Omega)}+C_1\intot\left(\left(\frac{t-s}{2}\right)^{-\alpha}+\left(\frac{t-s}{2}\right)^{-\frac{1}{2}-\alpha-(\frac{1}{p}-\frac{1}{q})}\right)e^{-\lambda_1\frac{t-s}{2}}\|u(\cdot,s)\nabla v(\cdot,s)\|_{\Lo^p(\Omega)}\intd s\nonumber
\end{align}
for $t\in(0,T_{max})$, where here and below every $C_i$ denotes a suitable positive constant not depending on $t$. To estimate $\|u\nabla v\|_{\Lo^p(\Omega)}$, we choose $r>1$ such that $pr<\theta$ still holds and $r'$ as corresponding H\"older conjugate. Then by H\"older's inequality we have
\begin{align}\label{propuinf-eq2}
\owncount
\|u\nabla v\|_{\Lo^p(\Omega)}\leq\|u\|_{\Lo^{pr'}(\Omega)}\|\nabla v\|_{\Lo^{pr}(\Omega)}\leq C_2\|\nabla v\|_{\Lo^\theta(\Omega)}\|u\|_{{\Lo^{pr'}(\Omega)}}\ \text{for all }t\in(0,T_{max}),
\end{align}
wherein by Proposition \ref{propnablavtheta} $\|\nabla v\|_{\Lo^\theta(\Omega)}$ is bounded independently of $t$ and thus
\begin{align*}
\|u\nabla v\|_{\Lo^p(\Omega)}\leq C_3\|u\|_{\Lo^{pr'}(\Omega)}\,\text{ for all }t\in(0,T_{max}).
\end{align*}
Further application of H\"older's inequality results in
\begin{align}\label{propuinf-eq3}
\owncount
\|u\|_{\Lo^{pr'}(\Omega)}\leq\|u\|_{\Lo^1(\Omega)}^{1-\beta}\|u\|_{\Lo^\infty(\Omega)}^{\beta},
\end{align}
where $\beta:=1-\frac{1}{pr'}\in(0,1)$. Next we define $M(T):=\sup_{t\in(0,T)}\|u(\cdot,t)\|_{\Lo^\infty(\Omega)}$, which for every $T\in(0,T_{max})$ is obviously finite, and recall $u_0\in C^0(\Omega)$ as well as $\|u\|_{\Lo^1(\Omega)}=\|u_0\|_{\Lo^1(\Omega)}$ (see Lemma \ref{conservmass}). Consequently \eqref{propuinf-eq1} -- \eqref{propuinf-eq3} imply
\begin{align*}
\|u(\cdot,t)\|_{\Lo^\infty(\Omega)}&\leq C_4+C_5\|u_0\|_{\Lo^1(\Omega)}^{1-\beta}M^\beta(T)\int\limits_0^\infty\left(\sigma^{-\alpha}+\sigma^{-\frac{1}{2}-\alpha-(\frac{1}{p}-\frac{1}{q})}\right)e^{-\lambda_1\sigma}\intd\sigma\quad\forall T\in(0,T_{max}).
\end{align*} Since we have $-\alpha>-1$ and $-\frac{1}{2}-\alpha-(\frac{1}{p}-\frac{1}{q})>-1$, the integral on the right hand side converges and it follows that
\begin{align*}
M(T)&\leq C_4+C_6 M^\beta(T)\quad\forall T\in(0,T_{max}).\end{align*}
Hence, since $\beta\in(0,1)$, by Lemma \ref{Malpha-ineqlem} we have $M(T)\leq C_7$ for all $T\in(0,T_{max}).$ This is equivalent to $\|u(\cdot,t)\|_\infty\leq C_7$ for all $t\in(0,T_{max})$ and thus completes the proof.
\end{bew}
In view of \eqref{eq2.1} the last two propositions ensure $T_{max}=\infty$ under our premise of $\sintomega u_0\intd x<4\pi$. Let us promptly gather these statements in the proof for the main theorem:
\begin{proof}[\textbf{Proof of Theorem \ref{globex2th}:}]
Combine propositions \ref{propnablavtheta} and \ref{propuinf} with the alternative \eqref{eq2.1} from Theorem \ref{localexist}.
\end{proof}

Recalling that in \cite{NSY97} it was shown that the bound for $m_u$ can be extended to $8\pi$, if $\Omega\subset\mathbb{R}^2$ is radial and the functions $u_0,v_0$ are radially symmetric, it seems sensible to assume a similar result holds in our case. This is indeed possible if we also require that $f$ is radially symmetric in order to ensure that the solutions to \eqref{KSf} remain symmetric as well. The only adjustment necessary in our proofs is taking a different -- but still similar -- inequality to the one stated in Theorem \ref{trudingermoser}. We briefly mentioned the inequality in question during the introduction to Section \ref{cglob2s1} and it can also be found in \cite[Theorem 2.1]{NSY97}. All remaining steps are left unchanged.

\setcounter{gleichung}{0} 

\section{Boundedness of solutions in the higher dimensional case}\label{cglobgen}
For the chemotaxis-system \eqref{KS} in the case $n\geq3$ it is known that, in contrast to the case $n=2$, a smallness condition of the mass $\sintomega u_0\intd x$ is not enough to ensure that all solutions to \eqref{KS} are globally bounded. On the other hand, if $\|u_0\|_{\Lo^p(\Omega)}$ and $\|\nabla v_0\|_{\Lo^\theta(\Omega)}$ are bounded by a sufficiently small number with $\theta>n$ and $p>\frac{n}{2}$, then every solution to \eqref{KS} remains bounded in time. Both of these results were presented by Winkler in \cite{win10jde}.

\subsection{Boundedness of solutions in higher dimensions~- existence and asymptotic properties}\label{cglobgens1}
\sectionmark{Boundedness of  solutions in higher dimensions}
A similar result to the one stated above holds if we also impose a smallness condition on $f$. We call a solution satisfying such a smallness condition with regard to $u_0,v_0$ and $f$ a small-data solution. Throughout this section we denote the local classical solution to \eqref{KSf} again by $(u,v)$ and let the parameters $n,\theta$ and $\delta_0$ satisfy the inequalities $2\leq n$, $0<\delta_0<1$ and $n<\theta<\frac{n^2+2n\delta_0}{n-2\delta_0}$. Additionally we assume $u_0\in C^0(\bomega)$, $v_0\in\W^{1,\theta}(\Omega)$, $f\in\Lo^\infty([0,\infty);\Lo^{\frac{n}{2}+\delta_0}(\Omega))\cap C^\alpha(\Omega\times(0,\infty))$, and that all necessary conditions of Theorem \ref{localexist} hold. In order to prove the boundedness of $\|u\|_{\Lo^\infty(\Omega)}$ and $\|v\|_{\W^{1,\theta}(\Omega)}$, we work along the same lines as in \cite{win10jde} using semigroup estimates and the following lemma (cf.\cite[Lemma 1.2]{win10jde}) to exclude the possibility of $T_{max}$ being finite and thus, in view of the alternative \eqref{eq2.1}, providing the bounds for the norms in question.
\begin{lem}\label{semigroupintegrallem}
Let $\alpha<1,\beta<1$ and $\gamma\neq\delta$ be positive constants. Then there exists $C>0$ such that
\begin{align*}
\intot\left(1+(t-s)^{-\alpha}\right)e^{-\gamma(t-s)}\cdot\!(1+s^{-\beta})e^{-\delta s}\intd s\leq C\left(1+t^{\min\lbrace0,1-\alpha-\beta\rbrace}\right)e^{-\min\lbrace\gamma,\delta\rbrace t}
\end{align*}
holds for all $t>0$.
\end{lem}

Without further preparations we can show the main result of this section.

\begin{proof}[\textbf{Proof of Theorem \ref{globalexist3+}:}]
Let $q_0:=\frac{n}{2}+\delta_0$. By the assumption imposed on $\theta$ we have $\frac{1}{\theta}>\frac{1}{q_0}-\frac{1}{n}$. 
We fix some small $\varepsilon_0>0$ -- to be specified later -- and assume that \eqref{globalexist3+eq1} holds for some $\varepsilon\in(0,\varepsilon_0)$. Now we define
\begin{align*}
T&:=\sup\bigg\lbrace T^*\in(0,T_{max}]\ \bigg\vert\,\|u(\cdot,t)-\etd u_0\|_{\Lo^p(\Omega)}\leq\varepsilon\left(1+t^{-\frac{n}{2}(\frac{1}{q_0}-\frac{1}{p})}\right)e^{-\frac{\lambda_1}{r}t}+\varepsilon\\
&\hspace*{160pt}\text{ for all }t\in(0,T^*)\text{ and each }p\in[\theta,\infty]\bigg\rbrace\leq\infty.
\end{align*}
Then $T$ is well-defined and positive with $T\leq T_{max}$, where $T_{max}$ denotes the maximal existence time corresponding to $u_0$ and $v_0$ obtained from Theorem \ref{localexist}. By the variation-of-constants formula it holds that
\begin{align*}
v(\cdot,t)-e^{\frac{t}{\tau}(\Delta-1)}v_0&=\frac{1}{\tau}\intot\etsdm\left(u(\cdot,s)-\overline{u_0}+f(\cdot,s)\right)\intd s+\frac{1}{\tau}\intot\etsdm\overline{u_0}\intd s,
\end{align*}
where $\overline{u_0}:=\frac{1}{\vert\Omega\vert}\sintomega u_0\intd x$. In what follows every $C_i$ will again denote a generic positive constant not depending on $t$. Applying $\nabla$ to both sides, the second integral disappears and we utilize semigroup estimates to obtain
\begin{align}\label{globalexist3+proof1}
\owncount
\ \left\|\nabla\left(v(\cdot,t)-e^{\frac{t}{\tau}(\Delta-1)}v_0\right)\right\|_{{\Lo^\theta(\Omega)}}\leq&\ C_1\intot\left(1+(t-s)^{-\frac{1}{2}}\right)e^{-(\lambda_1+1)(t-s)}\|u(\cdot,s)-\overline{u_0}\|_{{\Lo^\theta(\Omega)}}\intd s\nonumber\\
+&\ C_2\intot\left(1+(t-s)^{-\frac{1}{2}-\frac{n}{2}(\frac{1}{q_0}-\frac{1}{\theta})}\right)e^{-(\lambda_1+1)(t-s)}\|f(\cdot,s)\|_{{\Lo^{q_0}(\Omega)}}\intd s\nonumber\\ =:&\ C_1 I_1+C_2 I_2\,\text{ for all } t\in(0,T).
\end{align}
For $I_1$ we control $\|u(\cdot,s)-\overline{u_0}\|_{\Lo^p(\Omega)}$ from above, using the definition of $T$, semigroup estimates and the fact that $e^{-\lambda_1 s}\leq e^{-\frac{\lambda_1}{r}s}$ holds for $r>1$ and $s\geq0$, according to
\begin{align}\label{globalexist3+proof3}
\owncount
\|u(\cdot,s)-\overline{u_0}\|_{\Lo^p(\Omega)}\leq&\,\|u(\cdot,s)-e^{s\Delta} u_0\|_{\Lo^p(\Omega)}+\|e^{s\Delta}\left(u_0-\overline{u_0}\right)\|_{\Lo^p(\Omega)}\nonumber\\
\leq&\,\varepsilon\left(1+s^{-\frac{n}{2}(\frac{1}{q_0}-\frac{1}{p})}\right)e^{-\frac{\lambda_1}{r}s}+C_3\|u_0\|_{\Lo^{q_0}(\Omega)}\left(1+s^{-\frac{n}{2}(\frac{1}{q_0}-\frac{1}{p})}\right)e^{-\frac{\lambda_1}{r}s}+\varepsilon\nonumber\\
\leq&\,C_4\varepsilon\left(1+s^{-\frac{n}{2}(\frac{1}{q_0}-\frac{1}{p})}\right)e^{-\frac{\lambda_1}{r}s}+\varepsilon
\end{align}
for all $s\in(0,T)$ and every $p\in[\theta,\infty]$. Proceeding as in the proof of Lemma \ref{v estimate} we can estimate
\begin{align}\label{globalexist3+proof2}
\owncount
C_2I_2&\leq C_2\varepsilon\int_0^\infty\left(1+\sigma^{-\frac{1}{2}-\frac{n}{2}(\frac{1}{q_0}-\frac{1}{\theta})}\right)e^{-(\lambda_1+1)\sigma}\intd \sigma\leq\ C_5\varepsilon\,\text{ for all }  t\in(0,T),
\end{align}
where we used the fact that $\frac{1}{\theta}+\frac{1}{n}>\frac{1}{q_0}$ implies $\frac{n}{2\theta}+\frac{1}{2}>\frac{n}{2q_0}$. In a similar fashion, we see that
\begin{align}\label{globalexist3+proof2.5}
\owncount
C_1\varepsilon\intot\left(1+(t-s)^{-\frac{1}{2}}\right)e^{-(\lambda_1+1)(t-s)}\intd s\leq C_6\varepsilon\,\text{ for all } t\in(0,T).
\end{align}
So in fact, combining \eqref{globalexist3+proof1} -- \eqref{globalexist3+proof2.5} results in
\begin{align*}
\ &\left\|\nabla\left(v(\cdot,t)-e^{\frac{t}{\tau}(\Delta-1)}v_0\right)\right\|_{\Lo^\theta(\Omega)}\leq\\&\ C_{7}\varepsilon\intot\left(1+(t-s)^{-\frac{1}{2}}\right)e^{-(\lambda_1+1)(t-s)}\cdot\left(1+s^{-\frac{n}{2}(\frac{1}{q_0}-\frac{1}{\theta})}\right)e^{-\frac{\lambda_1}{r}s}\intd s+C_8\varepsilon\,\text{ for all } t\in(0,T).
\intertext{Further estimating the first term by means of Lemma \ref{semigroupintegrallem} implies}
\ &\left\|\nabla\left(v(\cdot,t)-e^{\frac{t}{\tau}(\Delta-1)}v_0\right)\right\|_{\Lo^\theta(\Omega)}\leq C_{9}\varepsilon\left(1+t^{\min\lbrace0,1-\frac{1}{2}-\frac{n}{2}(\frac{1}{q_0}-\frac{1}{\theta})\rbrace}\right)e^{-\min\lbrace\lambda_1+1,\frac{\lambda_1}{r}\rbrace t}+C_{8}\varepsilon\end{align*}
for all $t\in(0,T)$, wherein $\frac{n}{2}(\frac{1}{q_0}-\frac{1}{\theta})<\frac{1}{2}$. We therefore obtain the following bound for $\|\nabla v(\cdot,t)\|_{\Lo^\theta(\Omega)}$:
\begin{align}\label{globalexist3+proof4}
\owncount
\left\|\nabla v(\cdot,t)\right\|_{\Lo^\theta(\Omega)}&\leq\left\|\nabla\left(v(\cdot,t)-e^{\frac{t}{\tau}(\Delta-1)}v_0\right)\right\|_{\Lo^\theta(\Omega)}+\left\|\nabla e^{\frac{t}{\tau}(\Delta-1)}v_0\right\|_{\Lo^\theta(\Omega)}\leq C_{10}\varepsilon\quad\forall t\in(0,T).
\end{align}
With the estimate for $\|\nabla v\|_{\Lo^\theta(\Omega)}$ at hand, we are now able to apply similar steps to the variation-of-constants formula for $u(\cdot,t)$. In fact, for all $t\in(0,T)$ we have
\begin{align*}
\ &\ \left\|u(\cdot,t)-\etd u_0\right\|_{\Lo^p(\Omega)}\leq\intot\left\|\etsd\nabla\divdot(u(\cdot,s)\nabla v(\cdot,s))\right\|_{\Lo^p(\Omega)}\intd s\\
\leq&\ C_{11}\!\intot\!\left(1+(t-s)^{-\frac{1}{2}-\frac{n}{2}(\frac{1}{\theta}-\frac{1}{p})}\right)e^{-\lambda_1(t-s)}\left(\|u(\cdot,s)-\overline{u_0}\|_{\Lo^\infty(\Omega)}+\|\overline{u_0}\|_{\Lo^\infty(\Omega)}\right)\|\nabla v(\cdot,s)\|_{\Lo^\theta(\Omega)}\intd s.
\end{align*}
Now invoking the inequalities from \eqref{globalexist3+proof3} and \eqref{globalexist3+proof4} leads to
\begin{align*}
\left\|u(\cdot,t)-\etd u_0\right\|_{\Lo^p(\Omega)}
&\leq\ C_{12}\varepsilon^2\intot\left(1+(t-s)^{-\frac{1}{2}-\frac{n}{2}(\frac{1}{\theta}-\frac{1}{p})}\right)e^{-\lambda_1(t-s)}\left(1+s^{-\frac{n}{2q_0}}\right)e^{-\frac{\lambda_1}{r}s}\intd s\\
&+\ C_{13}\varepsilon^2\intot\left(1+(t-s)^{-\frac{1}{2}-\frac{n}{2}(\frac{1}{\theta}-\frac{1}{p})}\right)e^{-\lambda_1(t-s)}\intd s\,\text{ for all } t\in(0,T).
\end{align*}
Because of $\frac{1}{2}+\frac{n}{2}(\frac{1}{\theta}-\frac{1}{p})\leq\frac{1}{2}+\frac{n}{2\theta}<1$ and $\frac{n}{2q_0}=\frac{n}{n+2\delta_0}<1$ we can apply Lemma \ref{semigroupintegrallem} once again to the first integral as well as the estimation process used in \eqref{globalexist3+proof2} to the second integral, to obtain
\begin{align}
\left\|u(\cdot,t)-\etd u_0\right\|_{\Lo^p(\Omega)}&\leq C_{14}\varepsilon^2\left(1+t^{\min\lbrace0,\frac{1}{2}-\frac{n}{2\theta}-\frac{n}{2}(\frac{1}{q_0}-\frac{1}{p})\rbrace}\right)e^{-\frac{\lambda_1}{r}t}+C_{15}\varepsilon^2\,\text{ for all }  t\in(0,T).\nonumber
\intertext{Now using $\frac{1}{2}-\frac{n}{2\theta}\geq0$, a straightforward case analysis shows that irrespective of the sign of $\frac{1}{2}-\frac{n}{2\theta}-\frac{n}{2}(\frac{1}{q_0}-\frac{1}{p})$ the inequality $
1+t^{\min\lbrace0,\frac{1}{2}-\frac{n}{2\theta}-\frac{n}{2}(\frac{1}{q_0}-\frac{1}{p})\rbrace}\leq 2\left(1+t^{-\frac{n}{2}(\frac{1}{q_0}-\frac{1}{p})}\right)$ holds and thus}
\owncount\label{globalexeist3+proof5}
\left\|u(\cdot,t)-\etd u_0\right\|_{\Lo^p(\Omega)}&\leq 2C_{16}\varepsilon^2\left(1+t^{-\frac{n}{2}(\frac{1}{q_0}-\frac{1}{p})}\right)e^{-\frac{\lambda_1}{r}t}+C_{15}\varepsilon^2\,\text{ for all }  t\in(0,T).
\end{align}
If $\varepsilon_0$ fixed at the beginning was chosen small enough, i.e. that it satisfies $\varepsilon_0<\max\left\lbrace\frac{1}{2C_{16}},\frac{1}{C_{15}}\right\rbrace$, then the possibility of $T$ being finite is excluded by the continuity of $t\mapsto\left\|u(\cdot,t)-\etd u_0\right\|_{\Lo^p(\Omega)}$. Thus, the definition of $T$ immediately implies $T_{max}=\infty$. Therefore the solution $(u,v)$ to \eqref{KSf} exists globally in time and the alternative \eqref{eq2.1} entails the boundedness of the desired norms. Additionally \eqref{globalexist3+eq2} follows from \eqref{globalexeist3+proof5} applied to $p=\infty$. In order to verify \eqref{globalexist3+eq3} we apply the same steps as before, resulting in
\begin{align*}
\ &\ \left\Vert\nabla\!\left(v(\cdot,t)-e^{\frac{t}{\tau}(\Delta-1)} v_0-\frac{1}{\tau}\intot \etsdm e^{s\Delta}u_0\intd s\right)\right\Vert_{\Lo^\theta(\Omega)}\\
=&\ \left\Vert\frac{1}{\tau}\intot\nabla\etsdm\left(u(\cdot,s)- e^{s\Delta}u_0\right)\intd s+\frac{1}{\tau}\intot\nabla\etsdm f(\cdot,s)\intd s\right\Vert_{\Lo^\theta(\Omega)}\\
\leq&\ C_{17}\varepsilon^2\left(1+t^{\min\lbrace0,\frac{1}{2}-\frac{n}{2q_0}\rbrace}\right)e^{-\frac{\lambda_1}{r}t}+C_{18}(\varepsilon^2+\varepsilon)
\leq\ C_{19}\varepsilon^2 e^{-\frac{\lambda_1}{r}t}+C_{20}\varepsilon\,\text{ for all } t>1,
\end{align*}
which completes the proof.
\end{proof}

\end{document}